\documentclass[12pt]{article}

\usepackage{amsmath, amsthm, amssymb}

\newtheorem*{thm*}{Theorem}
\newtheorem{thm}{Theorem}
\newtheorem{cor}{Corollary}
\newtheorem{lem}{Lemma}

\newtheorem{prop}{Proposition}
\newcommand{\C}{\mathbb C}
\newcommand{\R}{\mathbb R}

\newcommand{\f}{\displaystyle \frac}

\theoremstyle{definition}

\theoremstyle{remark}

\DeclareSymbolFont{AMSb}{U}{msb}{m}{n}
\DeclareMathSymbol{\bbbn}{\mathalpha}{AMSb}{"4E}
\DeclareMathSymbol{\bbbz}{\mathalpha}{AMSb}{"5A}
\DeclareMathSymbol{\bbbr}{\mathalpha}{AMSb}{"52}
\DeclareMathSymbol{\bbbq}{\mathalpha}{AMSb}{"51}
\DeclareMathSymbol{\bbbc}{\mathalpha}{AMSb}{"43}

\newcommand{\Id}{\mathrm{Id}}

\numberwithin{equation}{section}

\title{On some exceptional cases in the integrability of the three--body problem }

\author{Alexei   V. Tsygvintsev \\ \\ \parbox{8cm}{\begin{center}
 Unit\'e de math\'ematiques pures et appliqu\'ees\\  
 Ecole Normale Sup\'erieure de Lyon\\ 46, all\'ee d'Italie, Lyon\\
F--69364  Lyon Cedex 07, France\\{\tt atsygvin@umpa.ens-lyon.fr}\end{center}}}

\begin{document}

\bibliographystyle{amsplain}

\maketitle

\begin{abstract}
\noindent 
We consider the Newtonian planar three--body problem with positive masses  $m_1$, $m_2$, $m_3$.
We prove that it does not have  an additional first integral meromorphic in the complex neighborhood of the parabolic Lagrangian orbit besides three exceptional cases $ \sum m_i m_j    / (\sum m_k)^2= 1/3$, 
$2^3/3^3$, $2/3^2$ where the linearized equations are shown to be partially integrable.
 This result completes the non--integrability analysis of the three--body problem started in  papers \cite{T6}, \cite{T1}  and based of the Morales--Ramis--Ziglin approach.

 \end{abstract}


\section{Introduction}

Let  $M\subset  \C^n$ be a complex domain and let $\mathrm{Hol}(M)$ be a set of functions $f \, :\, M\to \C$  holomorphic in $M$.
We consider a  system of $n$ ordinary differential equations written in the Pfaffiann form
\begin{equation} \label{foliation}
 \begin{array}{lll}
 \f{dx_1}{f_1(x)}=\f{dx_2}{f_2(x)}=\cdots=\f{dx_n}{f_n(x)}, \quad f_i\in \mathrm{Hol}(M)\,,\\ \\
x=(x_1 ,x_2,\dots,x_n)\in M.
\end{array}
\end{equation}
 
As far as the dynamical properties of the flow are concerned it is more convenient    to consider the trajectories of the vector field $F=(f_1,\dots,f_2)$  as leafs of a  codimension $n-1$  foliation defined  by \eqref{foliation} regardless  the time parametrisation.

Let $h \, :\, \Gamma \to M$  be  a  particular leaf of \eqref{foliation} where $h$   is a certain holomorphic map (not necessary unique).  We note that in many mechanical problems the surface  $\Gamma$  is a  punctured Riemann sphere $\Gamma=\C P^1\setminus  \{z_1,\dots,z_k \}$ and  $h$ is rational.  Let $e\in \Gamma$ be a fixed basepoint and $\{ \gamma\}$ is the set of loops generating   the fundamental group $\pi_1(e,\Gamma)$.

  We take a particular closed path $\gamma_0\in \{ \gamma\}$ defined by a continuos map $\gamma_0\, :\,  [0,1]  \mapsto \Gamma$, $\gamma(0)=\gamma(1)=e$.  Each tangent vector  $v\in T_e M$ can be transported along $\gamma_0$ following the neighborhood leafs of $h(\Gamma)$ back to the  tangent space $T_eM$.  One obtains a linear representation of $\pi_1(e,\Gamma)$ into $\mathrm{GL}(n,\C)$ called the {\it monodromy} group $G$.

This group  measures the complexity of ``enrolling''  of the  heighborhood   to $h(\Gamma)$ leafs  and usually contains  strong obstacles to the existence  of first integrals of \eqref{foliation} meromorphic  in $M$.   The following lemma of  Ziglin \cite{Zig} establishes  the link between the integrability of \eqref{foliation}  and rational  invariants of $G$.

\begin{lem}  \label{lz} Let $\Phi_1\dots \Phi_k$ be a set of  functionally independent  first 
integrals of the differential system \eqref{foliation} which are meromorphic  in $M$. Then the monodromy group $G\subset  \mathrm{GL}(\C,n)$ admits  $k$ functionally independent homogeneous rational   invariants $I_1,\dots,I_k$.
\end{lem}  

In many mechanical problems the previous  lemma allows to reduce the initial integrability problem  to the question from the theory of invariants of finitely generated linear groups $G=<G_1,\dots,G_k>\subset \mathrm{GL}(n,\C)$. Once it is shown  that $G$ has no more than $p$ rational invariants one concludes that  the system \eqref{foliation} can not have more than $p$ functionally independent  first integrals meromorphic in $M$.  If $G$ does have a rational invariant,  the higher variational approach has to be applied (see \cite{Morales}). 

\section{The planar three--body problem }

It is natural to view the monodromy generators as linear transformations obtained through the solving of  the {\it normal variational equation}  of \eqref{foliation} along the particular solution $h(\Gamma)$ (see for definition \cite{Zig}).     This equation describes the linearization of the system \eqref{foliation} around the particular orbit  after reduction of all already known first integrals.

We consider three mass points $m_1>0$, $m_2>0$, $m_3>0$ in the plane which attract each other
according to the Newtonian law.   Using the   Whittaker variables   \cite{W}  the corresponding equations of motion can be written as a Hamiltonian system with $3$ degrees of freedom 

\begin{equation} \label{EQ} \dot q_r=\displaystyle\frac{\partial
H}{\partial p_r},\quad \dot p_r=
-\displaystyle\frac{
\partial H}{\partial q_r}, \quad (r=1,2,3)\,,
\end{equation}  with the Hamiltonian function $$
\begin{array}{ll}
H=\displaystyle \frac{M_1}{2}\left\{p_1^2+\displaystyle
\frac{1}{q^2_1}P^2\right\}+\displaystyle \frac{M_2}{2}(p_2^2+
p_3^2)+\displaystyle \frac{1}{m_3}\left\{p_1p_2-\displaystyle
\frac{p_3}{q_1}P\right\} -\displaystyle \frac{ m_1m_3}{r_1}-
\displaystyle \frac{m_3m_2}{r_2}- \displaystyle \frac{
m_1m_2}{r_3}, \\ \\ P=p_3q_2-p_2q_3-k, M_1=m_3^{-1}+m_1^{-1},
M_2=m_3^{-1}+m_2^{-1},
\end{array}
$$  where $$ r_1=q_1, \quad r_2=\sqrt{q^2_2+q^2_3}, \quad
r_3=\sqrt{ (q_1-q_2)^2+q^2_3},$$ are the mutual distances of the
bodies; $q_1$ is the distance $m_3m_1$; $q_2$ and $q_3$ are the
projections of $m_2m_3$ on, and perpendicular to $m_1m_3$; $p_1$
is the component of momentum of $m_1$ along $m_3m_1$; $p_2$ and
$p_3$ are the components of momentum of $m_2$ parallel and
perpendicular to $m_3m_1$; $k$ is the constant of the angular
momentum.

 Let us denote $\mathrm{M}(4,K)$ the set of square $n$ by $n$ matrices over a  field $K$.
In \cite{T1} we calculated the normal variational equation   of \eqref{EQ} along  the Lagrangian parabolic equilateral solution for a fixed non-zero value of the angular momentum $k$

\begin{equation}\label{nve}  \displaystyle\frac{dx}{dz}=\left(
\displaystyle\frac{A}{z-z_0}+\displaystyle\frac{B}{z-z_1}+\displaystyle\frac{C}{z-z_2}\right)x, \quad x\in \C^4, \quad z \in \Gamma, \quad A,B,C\in M(4,\C)\,,
\end{equation}

where

\begin{equation}
 \begin{array}{lll}
 z_0=\displaystyle\frac{\sqrt{3}m_1m_2}{2S_2},\quad
z_1=\displaystyle\frac{m_1(\sqrt{3}m_2+iS_3)}{2S_2}, \quad
z_2=\displaystyle\frac{m_1(\sqrt{3}m_2-iS_3)}{2S_2},\\ \\
\Gamma=\C \setminus \{z_0,z_1,z_2 \}\\ \\
S_2=m_1m_2+m_2m_3+m_3m_1, \quad S_3=m_2+2m_3\,.
\end{array}
\end{equation}

 One verifies  (see \cite{T6}, Appendix A) that  in \eqref{nve}
 
\begin{equation}  \label{sym}
\left \{
\begin{array}{lll}
z_2=\overline z_1,\\ \\
A\in\mathrm{M}(4,\R), \quad B=R+iJ, \quad C=\overline B=R-iJ \quad  \mathrm{ with} \quad  R,J\in \mathrm{M}(4,\R)\,.
\end{array}
\right.
\end{equation}

Therefore,   the equations  \eqref{nve} are invariant under the complex conjugation fixing the time $t$. This is not surprising since the equations of the  three--body problem  \eqref{EQ} are real.

Let $\Sigma(z)$, $\Sigma(e)=\mathrm{Id}$, $e\in \Gamma$ be the  fundamental matrix solution of the linear differential equation \eqref{nve}.
 Continued along a closed path $\gamma\in \pi_1(e,\Gamma)$ the solution $\Sigma(z)$ gives a function $\tilde \Sigma(z)$ which also satisfies \eqref{nve}.  From linearity of \eqref{nve}  it follows that there
exists  $T_{\gamma}\in \mathrm{GL}(4,\C)$ such that $\tilde
\Sigma(z)=\Sigma(z) T_{\gamma}$. The set of matrices
$G=\{T_{\gamma}\}$ corresponding to all paths from  $\pi_1(e,\Gamma)$  form  monodromy group of the linear  system  \eqref{nve}. Let $T_i$ be the elements of $G$
corresponding to loops  around the singular points $z=z_i$,
$i=0,1,2$. Then $G$ is generated  by $T_0$, $T_1$,
$T_2$. Let $T_{\infty}\in G$  denotes the monodromy  element   around   $z=\infty$. 

\begin{prop}[\cite{T6}]
 The following assertions about the monodromy group $G$ hold

\noindent a)  The singularity $z_0$ is an apparent one i.e $T_0=\mathrm{Id}$ and 
\begin{equation} \label{Tinf}
T_1T_2=T^{-1}_{\infty}\,.
\end{equation}

\vspace{0.5cm}

\noindent b) The generators $T_1$, $T_2$ are unipotent trtansformations. Moreover, there exist  $U,V \in \mathrm{GL}(4,\C)$  such
that \begin{equation} \label{jordan} U^{-1}T_1U=V^{-1}T_2V=\left(
\begin{array}{cccc}
1 & 1 & 0 & 0 \\ 0 & 1 & 0 & 0 \\ 0 & 0 & 1 & 1
\\ 0 & 0 & 0 & 1
\end{array}
\right).
\end{equation}

\vspace{0.5cm}

\noindent c) 
\begin{equation}\label{TT}
{\rm{Spectr}}(T_{\infty})=\left\{ e^{2\pi i\lambda_1},\quad
e^{2\pi i\lambda_2},\quad e^{-2\pi i\lambda_1},\quad e^{-2\pi
i\lambda_2}\right\} \end{equation}  where \begin{equation}\label{lambdas}
\lambda_1=\displaystyle\frac{3}{2}+\displaystyle\frac{1}{2}\sqrt{13+\sqrt{\theta}},\quad
\lambda_2=\displaystyle\frac{3}{2}+\displaystyle\frac{1}{2}\sqrt{13-\sqrt{\theta}},\end{equation}
and \begin{equation} \theta=144\left(1-\displaystyle
3\frac{S_2}{S_1^2}\right).\end{equation} 
\end{prop}

  Let $U_{\Gamma}$ be a connected neighborhood of the Lagrangian parabolic equilateral solution of  the three--body problem \eqref{EQ}  defined in  \cite{T6}. 
   Below we summarize the results concerning the integrability of \eqref{EQ} obtained in our previous papers.

\begin{thm}{(\cite{CRAS1},\cite{T6})} \label{1}
For arbitrary  $m_i>0$, $i=1,2,3$ the monodromy group $G$ of \eqref{nve} does not have two independent rational invariants. As a result,   the  three--body problem  \eqref{EQ} never has two additional independent  first integrals meromorphic in $U_{\Gamma}$.
\end{thm}

\begin{cor}
The three--body problem \eqref{EQ}  is not completely meromorphically integrable in the sense of Liouville--Arnold.\end{cor}

We introduce the parameter  \begin{equation} \sigma=\f{m_1m_2+m_2m_3+m_3m_1}{( m_1+m_2+m_3)^2}\,. \end{equation}
The following two theorems  concern the non--existence of \emph{one} additional meromorphic first integral.

\begin{thm}{(\cite{T11},\cite{T1})} \label{2}
Let  $\sigma\not
\in \left\{ 
\frac{1}{3},  \frac{2^3}{3^3}, 
\frac{2}{9},  \frac{7}{48}, 
\frac{5}{2^4}\right\}.$
Then  the three--body problem \eqref{EQ} does not have an  additional
first integral meromorphic in $U_{\Gamma}$.\\
\end{thm}

\begin{thm}{(\cite{T11},\cite{T1})} \label{3}
 If $\sigma\in
\left\{\frac{1}{3},\frac{2^3}{3^3}\right\}$ then $G$ has a polynomial  invariant and 
the normal variational equation   \eqref{nve}  admits a first integral $I(x,z)$ which is
  a  quadratic polynomial  with respect to $x$ and which  is a rational function
with respect to $z$.
\end{thm}

The proof of these results is based on   the following result  from \cite{Zig}:  to every additional first integral of \eqref{EQ} independent with $H$ and meromorphic in $U_{\Gamma}$ corresponds a rational invariant of the monodromy group $G$.  We  used also  the infinitesimal  techniques from the Morales--Ramis differential Galois approach \cite{Morales}.

As follows  from Theorems \ref{1}--\ref{3}, the only remaining values of  $\sigma$ for those  the integrability property was not clear are 
\begin{equation} \label{main}
\sigma \in \left\{\frac{2}{9},  \frac{7}{48}, 
\frac{5}{2^4}\right\}\,.
\end{equation}

Our main result makes more precise the integrability property for these  values of $\sigma$.

\begin{thm}\label{MAIN}
 If  $\sigma=\frac{7}{48}, \frac{5}{2^4}$  then the three--body problem \eqref{EQ}  does not have an  additional first integral meromorphic in $U_{\Gamma}$. If $\sigma= \frac{2}{9}$ then $G$ has a polynomial invariant of degree $1$ or $2$ so that the normal variational equation \eqref{nve}  admits a first integral $I(x,z)$ which is
  a linear or   quadratic polynomial  with respect to $x$ and which  is a rational function
with respect to $z$.

\end{thm}

The proof is contained  in the next section.

\section{The reflection symmetry of the monodromy group}

As shown in \cite{T11}, \cite{T1} the monodromy group  $G$ always possesses a centralizer in $\mathrm{GL}(4,\C)$ given explicitly by 

\begin{equation} \label{c}
T=T_{\infty}+T_{\infty}^{-1}-2\, \mathrm{Id}\,. 
\end{equation}

\begin{prop}
Let 
\begin{equation} \label{r}
\sigma \in \left \{ \frac{7}{48}, \frac{5}{2^4}\right \}\,. \end{equation} 

Then $\mathrm{Spectr}(T)=\{ \sigma_1,\sigma_1, \sigma_2, \sigma_2\}$ where $\sigma_1 \neq \sigma_2$.
\end{prop}

It can be verified directly  with help of  the following formulas obtained  in \cite{T1} 

\begin{equation}
\mathrm{Spectr}(T)=\{ \sigma_1,\sigma_1, \sigma_2, \sigma_2\}, \quad \sigma_i=2(\cos(2\pi \lambda_i)-1), \quad i=1,2\,,
\end{equation}
where $\lambda_i$ are defined by \eqref{lambdas}.

Thus, if \eqref{r} holds, then the Jordan canonical form of $T$ always contains two $2\times 2$ blocks (either diagonal or not)  corresponding to two different eigenvalues $\sigma_1$ and $\sigma_2$.
Then, as follows from the solution of the Frobenius problem (see f.e \cite{G}),  the relations $[T_i,T]=0$, $i=1,2$ imply  the existence of a linear basis  in which the monodromy generators $T_{1,2}$ have  the same block--diagonal form
\begin{equation}
T_1=
\left[
\begin{array}{lll}
T^{11}_1&0\\
0&T^{22}_1
\end{array}
\right], \quad 
T_2=
\left[
\begin{array}{lll}
T^{11}_2&0\\
0&T^{22}_2
\end{array}
\right]\,, 
\end{equation} 
with unipotent blocks $T_i^{kk}\in \mathrm{GL}(2, \C)$.

\begin{prop} \label{prop1}
Under the condition  \eqref{r}  the group  $G$ does not have a  rational invariant. 
\end{prop}
\begin{proof}
The proof is exactly the same as in \cite{T6}, pp. 243-244.   Since $T\neq \alpha \mathrm{Id}$, $\alpha \in \C$ it is sufficient to verify that  $1\not \in \mathrm{Spectr}(T_{\infty})$ (*).   For  $\sigma=7/48$ one obtains $\lambda_1=3/2+\sqrt{22}/2$, $\lambda_2=5/2$ and for $\sigma=5/2^4$ respectively  $\lambda_1=7/2$, $\lambda_2=3/2+\sqrt{10}/2$.  In both cases $\lambda_i\not \in \mathbb Z$ and the condition (*)  follows from \eqref{TT}.
\end{proof}

We note  that the result of Proposition \ref{prop1} was not  contained in our work \cite{T6} where  $T_{\infty}$  was supposed  to be  diagonalizable.

\begin{lem}  \label {mainL} If $\sigma=\frac{2}{9}$ then $G$ has  a polynomial invariant.
\end{lem}

\begin{proof}

Below we assume  that the monodromy group $G=<T_1,T_2>$ of \eqref{nve}  is defined for the basepoint $e=0$.  One significant problem in the analysis of $G$ is  that the Jordan canonical form of $T_{\infty}$ and whose of the centralizer $T$ depend  on the masses $m_i$ in a quite complicated way. At least no elementary algebraic description of this dependence is known.  This difficulty can be overcome using the ``reflection'' symmetry of $G$ given by  the following lemma.

\begin{lem} \label{permut1}
The monodromy transformations $T_1$ and $T_2$ are related by $\overline T_1^{-1}=T_2$.
\end{lem}

\begin{proof}
 Let $\Sigma(z)$, $\Sigma(0)=\Id$ be the normalised fundamental matrix solution of \eqref{nve} and let $G$ be the corresponding monodromy group.  For a function $f(z)$ defined  by its Taylor expansion $f(z)=\sum  a_n z^n$, $a_n \in  \C$ we  define the operator  of complex conjugation $f(z)\stackrel{-}{\mapsto} \overline f(z)$ according to  $\overline f(z)=\sum \overline a_n z^n$.  One can always represent  locally $\Sigma(z)$ as a  power series convergent  in a small neighborhood of the regular point $z=0$. The symmetry conditions \eqref{sym}  imply then $\overline \Sigma(z)=\Sigma(z)$.   
  
  Let $\gamma_1$, $\gamma_2$  be two  loops  starting from  $z=0$ and going in the counter clock--wise direction around the singular points    $z_1 \in \C_+$ and  $z_2\in \C_-$ respectively.  By definition of $T_1$ we have (A): $\Sigma(z) \stackrel{\gamma_1}{\to} \Sigma(z) T_1$ after the analytic continuation of $\Sigma(z)$ along $\gamma_1$. Let $\overline \gamma_1$ denotes the loop symmetric to $\gamma_1$ with respect to the real axis $\mathrm{Im}\, z=0$ ( $\overline \gamma_1$  has the clock--wise orientation).  According to \eqref{sym} and (A) one will have   $\Sigma(z) \stackrel{\overline \gamma_1}{\to} \Sigma(z) \overline T_1$ after the ananlytic continuation along $\overline \gamma_1$.  At the same time we have (B): $\Sigma(z)  \stackrel{\gamma_2}{\to}  \Sigma(z)  T_2$  with $\Sigma(z)$  continued along the loop $\gamma_2$. The curves $\gamma_2$ and $\overline \gamma_1$ are homotopic   of  opposite  orientations. Hence, comparing (A) and (B) we obtain  $T_2=\overline T_1^{-1}$ that achieves  the proof of Lemma \ref{permut1}.
\end{proof}

If $\sigma=2/9$ then $\mathrm{Spectr}(T_{\infty})=(p,p,p^{-1},p^{-1})$, $p=e^{2\pi i\sqrt{3}}$, $p\neq p^{-1}$ in view of  \eqref{TT}. We denote by $L_p$ and $L_{p^{-1}}$ the eigenspaces of $T_{\infty}$ corresponding to the eigenvalues $p$ and $p^{-1}$ respectively. Firstly, we  consider the case  $\dim(L_p)=2$ (the case $\dim(L_{p^{-1}})=2$ is similar).  From Lemma \ref{permut1} and \eqref{Tinf} we deduce  $T_{\infty}=\overline T_{\infty}^{-1}$.  Therefore, if $v\in L_p$ then $\overline v\in L_{p^{-1}}$ and hence $\dim(L_p)=\dim(L_{p^{-1}})=2$. Thus, $T_{\infty}\in M(4,\C)$ is diagonalizable.   In this case, as shown in \cite{T1}, pp. 245-246, the group $G$ possesses a quadratic polynomial invariant that proves the statement. Let us consider the remaining case  $\dim(L_{p})=\dim(L_{p^{-1}})=1$ with $L_{p}=<v>$  and $L_{p^{-1}}=<\overline v>$ spanned  by the linearly independent vectors $v$ and  $\overline v$ respectively. Then the following proposition holds.

\begin{prop}
There exist  two linearly independent vectors $w$, $\overline w\in \C^4$ such that the   dual transformations $\mathcal{T}_1=T_1^T$, $\mathcal{T}_2=T_2^T$  act on $w$, $\overline w$ by permutations: 
\begin{equation}\label{rty}
\mathcal{T}_1\,  w=\overline{w}, \quad \mathcal{T}_1 \overline{w}= w, \quad  \mathcal{T}_2\, w=\overline{ w}, \quad \mathcal{T}_2 \overline{ w}= w\,.
\end{equation}
\end{prop}

\begin{proof}
Let $\mathcal{T}_{\infty}={T}_{\infty}^T$ then $\mathrm{Spectr}(\mathcal{T}_{\infty})=\mathrm{Spectr}({T}_{\infty})$. Since the geometric multiplicity of any eigenvalue in $\mathcal{T}_{\infty}$ is the same as in 
${T}_{\infty}$ we define $w$ and $\overline w$ as the only eigenvectors of  $\mathcal{T}_{\infty}$ corresponding to the eigenvalues $p$ and $p^{-1}$ respectively. The dual  centralizer  $\mathcal{T}=T^T$ commutes with $ \mathcal{T}_{1}$, $\mathcal{T}_{2}$ and its only eigendirections are $w$ and $\overline w$ as follows from  \eqref{c}.  The transformations $\mathcal{T}_{1}$, $\mathcal{T}_{2}$ preserve the eigendirections of $\mathcal{T}$.
So, in view of  \eqref{jordan}  and the identity $\mathcal{T}_{2}= 
\overline{\mathcal{T}_{1}}^{-1}$, either  the relations \eqref{rty} take place  or  we have
\begin{equation} \label{ret}
\mathcal{T}_{1} \, w= w,\quad  \mathcal{T}_{1} \, \overline w= \overline w,\quad  \mathcal{T}_{2} \, w= w,\quad  \mathcal{T}_{2} \, \overline w= \overline w\,. 
\end{equation}
 One defines  $<u,l>=\sum_{i=1}^n u_i l_i$ for  $u,l\in \C^4$. In the case \eqref{ret}  the monodromy group $G$ has two independent  polynomial invariants:  $I(x)=<w,x>$ and $I(x)=<\overline w,x>$. But it is impossible according to Theorem \ref{1}. The proposition  is proved. \end{proof}
 
In the remaining case  \eqref{rty}  the group $G$ has a linear  polynomial invariant given by $I(x)=<w+\overline w,x>$. That achieves the proof of Lemma \ref{mainL}. \end{proof}

 Theorem \ref{MAIN} follows  immediately  from  Proposition \ref{prop1}, Lemma \ref{mainL}, Lemma \ref{lz} and the fact  (see  \cite{T1}, p. 246) that   to every rational invariant of $G$ corresponds  a single--valued first integral of the linear differential system \eqref{nve}. We believe that these integrals, existing for $\sigma=1/3$, $2^3/3^3$, $2/3^2$ may contribute towards a better understanding of the dynamics of the three--body problem in the vicinity of parabolic Lagrangian orbits.

\providecommand{\bysame}{\leavevmode\hbox
to3em{\hrulefill}\thinspace}

\end{document}